# Analyzing the concept of super-efficiency in data envelopment analysis: A directional distance function approach

Mahmood Mehdiloozad[a], Israfil Roshdi[b, *]

**Abstract** Based on the framework of the directional distance function, we conduct a systematic analysis on the measurement of super-efficiency in order to achieve two main objectives. Our primary purpose is developing two generalized directional measures of super-efficiency that completely resolve the crucial infeasibility issue, commonly arisen in the traditional super-efficiency measures. The secondary goal is to demonstrate that our directional super-efficiency models encompass the conventional ones as special cases. The proposed measures are advantageous because they circumvent biases in super-efficiency estimation due to input and output slacks. They are general, and satisfy several desirable properties, such as always feasibility, monotonicity, unit independence, and translation invariance.

**Keywords:** Data envelopment analysis, Directional distance function, Super-efficiency, Infeasibility, Ranking·

---

[*] Corresponding Author

I. Roshdi

Department of Mathematics, Science and Research Branch, Islamic Azad University, Tehran, Iran

e-mail: i.roshdi@gmail.com

M. Mehdiloozad

Faculty of Mathematical Sciences and Computer, Kharazmi (Tarbiat Moallem) University, Tehran 15618, Iran

e-mail: m.mehdiloozad@gmail.com

## 1 Introduction

Since its introduction by Charnes et al. [1, 2], data envelopment analysis (DEA) has been established as a useful method for estimating the relative efficiency of a set of homogeneous decision making units (DMUs) with multiple inputs and multiple outputs. According to the DMUs' efficiency scores, DEA classifies the DMUs into two diverse *efficient* and *inefficient* groups. Unlike the inefficient DMUs, the efficient ones cannot be ranked based on their efficiencies because of having the same efficiency score of unity. It is not, however, reasonable to claim that the efficient DMUs have the same performance in actual practice. Now, the question arises how to rank the efficient DMUs? To address this question, different methods have been developed to achieve complete ranking of these DMUs. Adler et al. [3] provided a review of ranking methods, and drew the conclusion that "whilst each technique is useful in a specialist area, no one methodology can be prescribed … as the complete solution to the question of ranking".

Pioneered by Banker and Gifford [4] and Banker et al. [5], the super-efficiency (SE) implies the possible capability of a DMU in increasing its inputs and/or reducing its outputs without becoming inefficient [6]. The measurement of SE is important due to its use in a variety of applications, such as: (1) identifying outliers [4, 5, 7-9], (2) ranking the extreme-efficient DMUs [10], (3) measuring technology and productivity changes [11], (4) solving two-person ratio efficiency games [12], (5) identifying extreme-efficient DMUs [13], (6) analyzing sensitivity of efficiency classifications [14-19], (7) making acceptance decision rules [20], (8) calculating efficiency stability region [21], and (9) overcoming truncation problems in second-stage regressions intended to explain variation in efficiency [22]. Andersen and Petersen (AP hereafter) [10] introduced the super-efficiency as a ranking methodology to differentiate the performance of *extreme-efficient* DMUs. They proposed two CRS- (constant returns to scale) and VRS- (variable returns to scale) based AP models by making modifications on the CCR [1] and BCC [23] models. A major concern related to the AP models is the crucial infeasibility issue, which has been of considerable interest amongst researchers. The infeasibility of a SE model occurs when an efficient DMU under evaluation cannot reach the frontier formed by the rest of DMUs via increasing the inputs or decreasing the outputs, depending on the orientation of the model. Conditions for the feasibility in the CRS- and VRS-based AP models have been provided in [13, 17, 24-26]. Furthermore, as briefly reviewed in Section 2, a great deal of research efforts has been devoted to tackle the infeasibility problem. However, these efforts have not succeeded when the input-output data are non-negative, with two exceptions [27, 28].

Based on the directional distance function (DDF) of Chambers et al. [29, 30], Ray [31] developed a SE model (we call it Ray) to obtain a Nerlovian–Luenberger measure of SE by adjusting both input and output levels. Although the Ray model does not pose a similar infeasibility problem in the VRS-based AP model, it fails in two special situations. First, no feasible solution exists if the zero input value is present in a DMU under evaluation and all other DMUs in the reference set are positive-valued in that input. Second, this model may produce a technologically impossible projection point with negative coordinates. To addresses the infeasibility issues arisen in the Ray model, Chen et al. [6] recently developed a modified version of the Ray model (we call it M-Ray), by using another direction vector in the DDF. Although the M-Ray model resolves the infeasibility issue, its corresponding SE measure is not is not *complete* in the sense of [32], as it fails to account for the input and output slacks (i.e., the input savings and output surpluses). As a result, the SE scores obtained by the M-Ray are downwardly biased for super-efficient DMUs. In addition, the M-Ray measure has neither unit invariance nor translation invariance property, and thus is not able to deal with negative data.

In order to circumvent all the above-mentioned problems, we systematically analyze the SE measurement within the framework of DDF, and achieve two main objectives. As our primary purpose, we develop two generalized directional SE measures, i.e., Fractional GDSE and Linear GDSE, which are always feasible under generalized returns to scale (GRS, [33]) assumption. To accomplish this, first, we develop a super directional distance function (SDDF) using the DDF. Based on the SDDF, we propose a radial directional SE (RDSE) measure. Considering both the input and output slacks as additional sources of SE, we then extend the RDSE into our aimed measures. As the second objective, we demonstrate that the directional SE models encompass the conventional ones as special cases. This objective is fulfilled by specifying proper direction vectors to the directional SE models presented in this paper. It is shown that each of the Fractional and Linear GDSE super-efficiencies can be interpreted as a product of the means of the input and output super-efficiencies. The proposed measures are always feasible, and satisfy various desirable properties, such as, completeness [32], strong monotonicity [27, 34], unit invariance [32, 35], translation invariance (under the VRS) [32, 35-37], and positive affine transformation invariance [38]. The important property of translation invariance enables our measure to effectively deal with negative data.

The most significant practical advantage of our proposed measures is that due to the flexibility in computer programming, the computations of SE scores can be easily provided in standard DEA software, which will greatly help the decision makers. The other practical advantage of using our proposed measures is that they allow for the incorporation of decision makers' preferences into SE assessment. This leads to the generation of more realistic measures of SE by generating weights for individual inputs and outputs that are consistent with the objectives of decision makers. In some practical applications, simultaneous radial/non-radial (proportional/non-proportional) expansion and contraction of *all* inputs and/or *all* outputs may be impossible. In this sense, some inputs (outputs) are required to change radially while the remainders changes non-radially. To deal with such situations, we formulate a hybrid directional SE measure, i.e., Linear HDSE, which is always feasible. Furthermore, some applications allow only for the modifications in inputs (outputs) with unchanged outputs (inputs). In order to deal with such situations, we formulate two radial and non-radial input-oriented directional SE measures, and provide necessary and sufficient conditions for their feasibility. Since the exposition of output-oriented models is identical to the input-oriented models, we do not discuss about the output-oriented models in this research.

The remainder of the paper unfolds as follows. Section 2 gives a brief review of the previous research efforts on the super-efficiency measurement. Section 3 deals with the description of technology, followed by an introduction of the SDDF. Section 4 first introduces a radial directional super-efficiency model; second, develops the Fractional and Linear GDSE measures; third, provides a detailed discussion about properties of the proposed measure; fourth, extends the Linear GDSE measure into the Linear HDSE measure; and, finally develops two input-oriented directional measures of super-efficiency. In Sections 5, we demonstrate that the directional models include the conventional ones as special cases. Section 6 presents three illustrative examples to explain the properties and advantages of the proposed measures over the existing procedures. Section 7 concludes with some remarks.

## 2 Previous Studies on the SE

### 2.1 CRS-based SE

The previous studies on the SE measurement under the CRS assumption are classified into two groups:

**(I) Radial:** The basic CRS-based SE model is the CRS-based AP model proposed in the seminal work of AP. In addition to the issue of infeasibility, the CRS-based AP model may be unstable, as it is extremely sensitive to small variations in the data when some DMUs have rather small values for some inputs. In order to overcome these drawbacks, Mehrabian et al. [39] introduced the MAJ model that is stable, but not always feasible. Saati et al. [40] extended the MAJ model into a non-oriented version (we call it R-MAJ), and asserted that their model is always feasible. However, as shown in Remark 4.1.1 and Example 1, the R-MAJ model may produce a technologically impossible projection point with negative coordinates. Jahanshahloo et al. [41] changed the technique used for making the MAJ model unit-invariant, and proposed a slightly modified version of this model (we call it M-MAJ).

**(II) Non-Radial:** A common drawback of the above-mentioned radial models is that none of them considers the possible non-zero slacks. To circumvent this problem, Tone [34] assessed the SE based on the SBM model [42], and introduced the so-called Super-SBM-I, Super-SBM-O and Super-SBM-C measures. The Super-SBM-I and Super-SBM-O models are not always feasible; the Super-SBM-C model requires positive input-output data for efficient DMUs. Using the L1-norm, Jahanshahloo et al. [28] presented a SE model (we call it Norm1), which is always feasible. Li et al. [43] developed the LJK model, and provided a certain condition for its feasibility. Based upon the additive DEA model [44], Du et al. [27] presented two additive (i.e., slack-based) SE models (we call them Super-Add (I) and Super-Add (II)). While the Super-Add (I) model is feasible for non-negative data, the original SE index associated with this model cannot be defined for efficient DMUs with some zero input and/or output values. To resolve this problem, they suggested two alternative indices. Thought the original Super-Add (II) model may not be feasible for non-negative data, two modified versions of it with alternative objective functions (we call them Super-Add (III) and Super-Add (IV)) are always feasible. As we show, the Super-Add (III) is equivalent to the Norm1 model and vice versa, and the only difference between these models lies in their associated indices. Based upon the ERM [45] model, Ashrafi et al. [46] developed a SE model (we call it Super-ERM), which is equivalent to the Super-SBM-C model and vice versa.

### 2.2 VRS-Based Super-Efficiency

The basic VRS-based SE model is the VRS-based AP model. Extending the Super-SBM-C model for the case of VRS, Tone [34] developed the Super-SBM-V model. Then, Lovell and Rouse [47] suggested using a user-defined scaling factor to make the VRS-based AP model feasible. Although their suggested approach provides feasible solutions for all efficient DMUs, it does not fully resolve the infeasibility problem since

- The provided SE scores for those efficient DMUs without feasible solutions in the VRS-based AP model (i.e., problematic DMUs) are equal to the user-defined scaling factor.
- Thus, these scores subjectively depend on the scaling factor, and have no economic interpretations.
- The problematic DMUs are all tied at the maximum score.

Jahanshahloo et al. [28] extended the Norm1 model to a VRS-based SE model (we call it Norm1-V). Chen [48, 49] developed a revised SE model to find out the possible SE in terms of input saving or output surplus for infeasibility by projecting the inefficient DMUs onto the VRS frontier via the BCC model. However, this revised model may still encounter infeasibility, as it based on the input- and output-oriented VRS-based AP super-efficiencies. Another deficiency of this approach, and that of Lovell and Rouse [47], is using the perturbated data for estimating SE scores. Johnson and McGinnis [8] developed a hyperbolic-oriented SE model, which is free from the infeasibility problem in the case of positive input-output data. However, this model is feasible only for some certain cases when the positivity assumption is relaxed. In addition to this drawback, the hyperbolic model is a nonlinear programming problem, and does not consider the non-zero possible sacks. Cook et al. [50] developed a two-stage process to address the infeasibility issue, and defined a SE score with respect to both input and output super-efficiencies. As we show, the feasibility of this approach cannot be guaranteed for non-negative data. Based on the simultaneous proportionate increase in inputs and decrease in outputs, Chen et al. [51] proposed a SE model, which is feasible for the positive data. Extending the work by Chen [49], Lee et al. [52] developed an alternative two-stage process for overcoming the infeasibility issue, which has been integrated into a single model by Chen and Liang [53].

## 3 Preliminaries

Let us start by introducing the necessary notations and defining the basic concepts used throughout the paper. We assume that we are dealing with $n$ observed DMUs; each uses $m$ inputs to produce $s$ outputs. We denote by $x_j = \left(x_{1j}, \dots, x_{mj}\right)^T \in \mathrm{IR}^m_{\geq 0}$ and $y_j = \left(y_{1j}, \dots, y_{sj}\right)^T \in \mathrm{IR}^s_{\geq 0}$ the input and output vectors of $\mathrm{DMU}_j, j \in \{1, \dots, n\}$, where $x_j \neq 0$ and $y_j \neq 0$. We consider $\mathrm{DMU}_o$ as the DMU under evaluation.

### 3.1 Technology

The production possibility set or technology $T$ is the set of all technologically feasible input-output vectors:

$$T = \{(x, y) | x \in \mathrm{IR}^m_{\geq 0} \text{ can produce } y \in \mathrm{IR}^s_{\geq 0}\}. \tag{1}$$

We assume that the technology satisfies the following assumptions: (T.1) no free lunch; (T.2) boundedness; (T.3) closedness; (T.4) free disposability of inputs and outputs; and (T.5) convexity (see [11] for details).

By excluding $\mathrm{DMU}_o$ from the reference set, the explicit DEA-based representation of the technology $T$ under the GRS assumption is as follows:

$$T = \left\{ (x, y) \in \mathrm{IR}^{m+s}_{\geq 0} \,\middle|\, x \geq \sum_{j \in J} \lambda_j x_j \,, y \leq \sum_{j \in J} \lambda_j y_j \,,\ L \leq \sum_{j \in J} \lambda_j y_j \leq U, \lambda_j \geq 0, j \in J \right\}, \tag{2}$$

where $J = \{1, \dots, n\} \backslash \{o\}$, and $L$ $(0 \leq L \leq 1)$ and $U$ $(1 \leq U)$ are upper and lower bounds for the sum of $\lambda_j$.

Note that the special cases $L = 0, U = \infty$ and $L = U = 1$ in (2) respectively correspond to the CRS- and VRS-based technologies.

### 3.2 Super Directional Distance Function

The technology can be completely characterized by the use of distance functions. An instance of such a function is the DDF [29, 30] that includes the Shephard [54, 55] distance functions as special cases. While the DDF seeks to simultaneously minimize the inputs and maximize the outputs of a given DMU in a pre-specified direction, the SE measurement requires the inputs to be maximized and the outputs to be minimized. We, therefore, make a slight modification on the DDF, and introduce the concept of super directional distance function (SDDF).

**Definition 3.2.1** Let $(x, y) \in \mathrm{IR}^m_{\geq 0} \times \mathrm{IR}^s_{\geq 0}$ be an input-output vector and $g = (g^-, -g^+) \in \mathrm{IR}^m_+ \times (-\mathrm{IR}^s_+)$ be a direction vector; the function

$$\vec{S}_g(x, y) := Min\{1 + \tau | (x + \tau g^-, y - \tau g^+) \in T\}. \tag{3}$$

is called *super directional distance function*.

In this definition, $g$ is a predetermined direction along which the vector $(x, y)$ is projected onto the frontier of $T$ at $\left(x + \vec{S}_g g^-, y - \vec{S}_g g^+\right)$; and $\vec{S}_g = \vec{S}_g(x, y)$ represents the minimum distance needed to reach the projection point in the direction $g$.

The SDDF seeks to measure the maximum input savings and output surpluses, which can be achieved by $(x, y)$ compared to the $T$. It has several useful properties, as given below.

*S.1.* $\vec{S}_g(x, y) = 1 - \vec{D}_g(x, y)$ where $\vec{D}_g(\cdot)$ is the DDF.

*S.2.* T-indication: $\vec{S}_g(x, y) \leq 1$ iff $(x, y) \in T$.

*S.3.* Translation: $\vec{S}_g(x + \alpha g^-, y - \alpha g^+) = \vec{S}_g(x, y) - \alpha$, for all $\alpha \in \mathrm{IR}$.

*S.4.* Homogeneity of minus one: $\vec{S}_{\lambda g}(x, y) = (1/\lambda)\vec{S}_g(x, y)$, for all $\lambda > 0$.

*S.5.* Input Monotonicity: $x \leq x'$ iff $\vec{S}_g(x, y) \geq \vec{S}_g(x', y)$.

*S.6.* Output Monotonicity: $y \leq y'$ iff $\vec{S}_g(x, y) \leq \vec{S}_g(x, y')$.

Property S.1 directly follows from the definitions of the DDF and SDDF. Furthermore, Properties S.2–P.6 are immediate consequences of property S.1 and the corresponding properties of the DDF.

# 4 Our Proposed SE Measures

## 4.1 Radial Directional SE Measure

For any given direction vector $g$, the $\vec{S}_g(\cdot)$ can be computed for $\mathrm{DMU}_o$ by solving the following problem, referred to as the *radial directional super-efficiency* (RDSE) model:

$$
\begin{aligned}
\gamma_o(g) := Min\ \ & 1+\tau \\
s.t.\ \ & \sum_{j\in J}\lambda_j x_{ij} \le x_{io} + \tau g_i^-, \quad i=1,\dots,m, \\
& \sum_{j\in J}\lambda_j y_{rj} \ge y_{ro} - \tau g_r^+, \quad r=1,\dots,s, \\
& L \le \sum_{j\in J}\lambda_j \le U, \lambda_j \ge 0, \forall j, \\
& y_{ro} - \tau g_r^+ \ge 0, \forall r, \tau{:}\, free\ in\ sign.
\end{aligned}
\tag{4}
$$

The optimal value of model (4), $\gamma_o(g)$, is referred to as the SE score of $\mathrm{DMU}_o$ if it is greater than the unity, and the bigger the value is, the more super-efficient it is. If $\mathrm{DMU}_o$ be *super-efficient*, i.e., $\gamma_o(g) > 1$, then the inputs and outputs of $\mathrm{DMU}_o$ have to be, respectively, expanded and contracted to get the projection $(x + \gamma_o g^-, y - \gamma_o g^+)$ in the technology $T$.

**Remark 4.1.1** *For having a "true" projection onto the $T$, i.e., $(x_o + \gamma_o g^-, y_o - \gamma_o g^+) \in T$, the constraint $y_{ro} - \tau g_r^+ \ge 0, r = 1, \dots, s$, should be added to the set of constraints of model (4).*

Let $T$ be a technology satisfying T1–T5, and $g = (g^-, -g^+) \in \mathrm{IR}_+^m \times (-\mathrm{IR}_+^s)$ be a chosen direction vector. *If $g^+ \neq 0$ and $s \ge 2$, then there exists some input-output vector $(x, y) \in IR_{\ge 0}^{m+s}$ such that the direction $g$ is infeasible at $(x, y)$, i.e., $\{(x, y) + \beta g | \beta \in IR\} \cap T = \emptyset$* (see [56] for details). From this result, we deduce that if $s \ge 2$ and the direction $g$ is infeasible at $\mathrm{DMU}_o$, then model (4) will be infeasible. This is illustrated via the following example.

**Table 1** The data set for example

| | $\mathbf{DMU_1}$ | $\mathbf{DMU_2}$ | $\mathbf{DMU_3}$ |
|---|---|---|---|
| $\mathbf{I_1}$ | 1 | 1 | 1 |
| $\mathbf{O_1}$ | 1 | 0 | 0 |
| $\mathbf{O_2}$ | 0 | 1 | 2 |

**Example 4.1.1** Table 1 shows three DMUs with one input and two outputs. Let us $g = (-g^-, g^+) = (-1,1,2)$ is the chosen direction vector and $\mathrm{DMU}_3$ is the under assessment unit. Then, the feasible region of the CRS case of the model (4) can be expressed as the following system of inequalities:

$$
\begin{aligned}
& \lambda_1 + \lambda_2 \le 1 + \tau \\
& \lambda_1 \qquad\ \ \le 0 - \tau \\
& \qquad\ \lambda_2 \ge 2 - 2\tau \\
& \lambda_1 \ge 0, \lambda_2 \ge 0, \\
& 0 - \tau \ge 0, 2 - 2\tau \ge 0.
\end{aligned}
\tag{5}
$$

Obviously, this system has no feasible solution, indicating that the CRS case of the model (4) is infeasible. However, removing the last two constraints (the additional constraints mentioned in Remark 4.1.1) in (5) leads to a consistent system, in which minimizing the objective function $1+\tau$ on the feasible region of (5) gives the optimal value of 1.3333 and the respective projection point $\text{Proj}(\text{DMU}_3) = \left(1.3333, \underline{-0.3333}, 1.3333\right)$, which is technologically impossible. This provides a theoretical justification for Remark 4.1.1.

As shown in Section 5, the R-MAJ and Ray models are two special cases of the model (4) without the constraints $y_{ro} - \tau g_r^+ \geq 0, r = 1, \dots, s$. In Examples 1 and 2, we demonstrate that neglecting theses constrains may lead to some feasible solutions where their corresponding projections, having some negative components, are technologically impossible.

### 4.2 Fractional and Linear Generalized Directional Super-Efficiency Models

As illustrated in the preceding subsection, the RDSE model in (4) may encounter infeasibility in some cases. Furthermore, the RDSE measure is not *complete* in the sense of [32], as it fails to account for the additional sources of SE caused by the non-zero slacks. As a result, RDSE model underestimates the SE of those extreme-efficient DMUs suffering from these slacks.

To circumvent these problems, we develop two generalized versions of RDSE model; each of them, being always feasible, reflects all the non-zero slacks in a single SE score. In doing so, we use the variables $\tau_i^-$, $i = 1, \dots, m$ and $\tau_r^+$, $r = 1, \dots, s$ as the individual rates of expansion and contraction in the $i$th input and the $r$th output of $\text{DMU}_o$, and introduce Fractional and Linear GDSE models as follows.

**Fractional GDSE**

$$\begin{aligned} \rho_o^F(g) \coloneqq Min \quad & \frac{1 + \frac{1}{m}\sum_{i=1}^{m} \tau_i^-}{1 - \frac{1}{s}\sum_{r=1}^{s} \tau_r^+} \\ s.t. \quad & \sum_{j \in J} \lambda_j x_{ij} \leq x_{io} + \tau_i^- g_i^-, \quad i = 1, \dots, m, \\ & \sum_{j \in J} \lambda_j y_{rj} \geq y_{ro} - \tau_r^+ g_r^+, \quad r = 1, \dots, s, \\ & L \leq \sum_{j \in J} \lambda_j \leq U, \\ & \lambda_j \geq 0, \tau_i^- \geq 0, \tau_r^+ \geq 0, \forall j, i, r. \end{aligned} \tag{6}$$

**Linear GDSE**

$$\begin{aligned} \varphi_L(g) \coloneqq Min \quad & \frac{1}{m}\sum_{i=1}^{m} \tau_i^- + \frac{1}{s}\sum_{r=1}^{s} \tau_r^+ \\ s.t. \quad & \sum_{j \in J} \lambda_j x_{ij} \leq x_{io} + \tau_i^- g_i^-, \quad i = 1, \dots, m, \\ & \sum_{j \in J} \lambda_j y_{rj} \geq y_{ro} - \tau_r^+ g_r^+, \quad r = 1, \dots, s, \\ & L \leq \sum_{j \in J} \lambda_j \leq U, \\ & \lambda_j \geq 0, \tau_i^- \geq 0, \tau_r^+ \geq 0, \forall j, i, r. \end{aligned} \tag{7}$$

Note that models (6) and (7) have the same constraints, but different objective functions. For a given direction vector $g$ , the objective functions of these models seek the minimum non-radial expansions in inputs and contractions in outputs, which are needed to reach the frontier constructed by the remaining DMUs. While Linear NDSE model is a linear programming problem, Fractional GDSE model is a non-linear programming problem. Nevertheless, one can transform the latter into an equivalent linear programming problem by using the "Charnes-Cooper transformation" [57].

Corresponding to Fractional and Linear GDSE models, we now introduce two SE measures, i.e., $\rho_o^F(g)$ and $\rho_o^L(g)$, where the former is given in the objective of (6) and the latter is defined as follows:

$$\rho_o^L(g) \coloneqq \left[1 + \frac{1}{m}\sum_{i=1}^{m} \tau_i^{-*}\right] \times \left[1 - \frac{1}{s}\sum_{r=1}^{s} \tau_r^{+*}\right]^{-1},$$

in which $\left(\lambda_j^*, \tau_i^{-*}, \tau_r^{+*}, \forall j, i, r\right)$ is an optimal solution to model (7). Since the constraints of models (6) and (7) are the same, the relationship $\rho_o^F(g) \leq \rho_o^L(g)$ holds.

An important feature of the measures $\rho_o^F(g)$ and $\rho_o^L(g)$ is circumventing biases due to input and output slacks in the SE estimation. Therefore, they generate accurate SE scores, and have a higher degree of discriminatory power than the RDSE measure.

Given an optimal solution $\left(\lambda_j^*, \tau_i^{-*}, \tau_r^{+*}, \forall j, i, r\right)$ to model (6) or (7), the projection $(\hat{x}, \hat{y})$ for $\text{DMU}_o$ is defined as

$$\begin{aligned} x_{io} &\to \hat{x}_{io} = x_{io} + \tau_i^{-*} g_i^-, \quad i = 1, \ldots, m, \\ y_{ro} &\to \hat{y}_{ro} = y_{ro} - \tau_r^{+*} g_r^+, \quad r = 1, \ldots, s. \end{aligned} \tag{8}$$

*Remark 4.2.1* The outputs of the projection points obtained from Fractional and Linear GDSE models are non-negative. Therefore, there is no necessity for adding the constraints $y_{ro} - \tau_r^+ g_r^+, r = 1, \ldots, s$ in these models (see Remark 4.1.1.).

*Remark 4.2.2* In models (6) and (7), the input and output inequality constraints corresponding to positive adjustment factors would be active at the optimum. In the other words, if $\tau_i^{-*} > 0$, $i \in \{1, \ldots, m\}$ or $\tau_r^{+*} > 0$, $r \in \{1, \ldots, s\}$, then the $i$th input or $r$th output constraints is binding.

*Remark 4.2.3* Positivity of the direction sub-vectors $g^-$ and $g^+$, i.e. $g^- > 0$ and $g^+ > 0$, is a sufficient condition for the feasibility of Fractional and Linear GDSE models. However, a necessary condition for the feasibility of these models to be imposed upon the chosen direction vector $g$ is that $g_r^+ > 0$ for any $r \in Q_o \coloneqq \left\{r \middle| \sum_{j \in J} y_{rj} = 0, y_{ro} > 0\right\}$ and that $g_i^- > 0$ for any $i \in P_o \coloneqq \left\{i \middle| x_{io} = 0, x_{ij} > 0, j \in J\right\}$.

*4.2.1 Well-definedness of Fractional and Linear GDSE measures:*

Given an optimal solution $(\lambda_j^*, \tau_i^{-*}, \tau_r^{+*}, \forall j, i, r)$ to model (6) (model (7)), the SE measure $\rho_o^F(g)$ ($\rho_o^L(g)$) would be well-defined, i.e. $\rho_o^F(g) \geq 1$ ($\rho_o^L(g) \geq 1$), provided that $\tau_r^{+*} \leq 1$ for all $r = 1, \dots, s$. In order to fulfill these conditions, consider the $r$th coordinate of $(\hat{x}, \hat{y})$ defined in (8). From Remark 4.2.1, we have $y_{ro} - \tau_r^{+*} g_r^+ \geq 0$, or equivalently $\tau_r^{+*} \leq \frac{y_{ro}}{g_r^+}$. Now, imposing the following condition on the chosen direction vector $g$ ensures the well-definedness of $\rho_o^F(g)$ ($\rho_o^L(g)$) under the GRS assumption:

$$\underset{r=1,\dots,s}{\text{Max}} \left\{ \frac{y_{ro}}{g_r^+} \right\} \leq 1. \tag{9}$$

For example, the condition (9) holds for the following direction vectors:

$$g_i^- = x_{io}, g_r^+ = y_{ro}, i = 1, \dots, m, r = 1, \dots, s; \tag{10}$$

$$g_i^- = \underset{j \in J}{\text{Max}} \{x_{ij}\}, g_r^+ = \underset{j \in J}{\text{Max}} \{y_{rj}\}, i = 1, \dots, m, r = 1, \dots, s. \tag{11}$$

The conditions (9) is sufficient for the well-definedness of the SE measure $\rho_o^F(g)$ ($\rho_o^L(g)$) under arbitrary RTS. However, one can impose another condition for assuring the well-definedness under a specific RTS. For example, in the VRS case, since $\sum_{j \in J} \lambda_j^* = 1$, the $r$th output constraint of model (6) results the following relationship:

$$\sum_{j \in J} \lambda_j^* y_{rj} \leq y_{ro} - \tau_r^{+*} g_r^+ \quad \Longrightarrow \quad \tau_r^{+*} \leq \frac{y_{ro} - \sum_{j \in J} \lambda_j^* y_{rj}}{g_r^+} \leq \frac{y_{ro} - \underset{j \in J}{\text{Min}} \{y_{rj}\}}{g_r^+}. \tag{12}$$

Based on (12), imposing the following condition on the $g$ ensures the well-definedness of the $\rho_o^F(g)$ ($\rho_o^L(g)$):

$$\underset{r=1,\dots,s}{\text{Max}} \left\{ \frac{y_{ro} - \underset{j \in J}{\text{Min}} \{y_{rj}\}}{g_r^+} \right\} \leq 1. \tag{13}$$

For example, as well as (10) and (11), the following direction vector fulfills the condition (13).

$$g_i^- = \underset{j \in J}{\text{Max}} \{x_{ij}\} - \underset{j \in J}{\text{Min}} \{x_{ij}\}, g_r^+ = \underset{j \in J}{\text{Max}} \{y_{rj}\} - \underset{j \in J}{\text{Min}} \{y_{rj}\}, i = 1, \dots, m, r = 1, \dots, s. \tag{14}$$

*4.2.2 Interpretation of Fractional and Linear GDSE measures:*

To better interpret the GDSE of $\text{DMU}_o$, having an optimal solution $(\lambda_j^*, \tau_i^{-*}, \tau_r^{+*}, \forall j, i, r)$ to model (6), the $\rho_o^F(g)$ can be expressed as

$$\rho_o^F(g) = \left[1 + \frac{1}{m}\sum_{i=1}^{m}\tau_i^{-*}\right] \times \left[1 - \frac{1}{s}\sum_{r=1}^{s}\tau_r^{+*}\right]^{-1}$$

$$= \left[\frac{1}{m}\sum_{i=1}^{m}\frac{g_i^- + \left((x_{io} + \tau_i^{-*}g_i^-) - x_{io}\right)}{g_i^-}\right] \times \left[\frac{1}{s}\sum_{r=1}^{s}\frac{g_r^+ - \left(y_{ro} - (y_{ro} - \tau_r^{+*}g_r^+)\right)}{g_r^+}\right]^{-1} \quad (15)$$

$$= \left[\frac{1}{m}\sum_{i=1}^{m}\frac{g_i^- + s_i^{-*}}{g_i^-}\right] \times \left[\frac{1}{s}\sum_{r=1}^{s}\frac{g_r^+ - s_r^{+*}}{g_r^+}\right]^{-1},$$

where $s_i^{-*} = (x_{io} + \tau_i^{-*}g_i^-) - x_{io}$, $i = 1, \dots, m$ and $s_r^{+*} = y_{ro} - (y_{ro} - \tau_r^{+*}g_r^+)$, $r = 1, \dots, s$, respectively represent input saving and output surplus of the $i$th input and the $r$th output in the direction $g$.

In the first term in the right-hand side, the ratio $\frac{g_i^- + s_i^{-*}}{g_i^-}$ represents the relative expansion rate of the $i$th input in the direction of $g_i^-$; hence, the first term is the mean expansion rate of inputs, and measures the *average input GDSE*. Similarly, in the second term, the ratio $\frac{g_r^+ - s_r^{+*}}{g_r^+}$ represents the relative contraction rate of the $r$th output in the direction of $g_r^+$. Thus, the second term is the mean contraction rate of outputs, and measures the *average output GDSE*.

In summary, the $\rho_o^F(g)$ can be well interpreted as the product of average input and output GDSEs, respectively given by $\rho_o^{FI}(g) \coloneqq 1 + \frac{1}{m}\sum_{i=1}^{m}\tau_i^{-*}$ and $\rho_o^{FO}(g) := \left[1 - \frac{1}{s}\sum_{r=1}^{s}\tau_r^{+*}\right]^{-1}$. A similar interpretation holds for the $\rho_o^L(g)$.

### 4.3 Properties of Fractional and Linear GDSE measures:

As mentioned earlier, Fractional and Linear NDSE models have several useful and interesting properties. By using suitable direction vectors, these models will have additional advantages as summarized below.

*Property1:* Always feasibility

**Theorem 4.3.1** *Models (6) and (7) are always feasible.*

***Proof*** Let $c \in [L, U]$ is a constant value, then we can define $\lambda_j' = b = \frac{c}{n-1}$, $j \in J$, $\tau_i^{-\prime} = Max\left\{0, \frac{1}{g_i^-}\left(-x_{io} + b\sum_{j\in J}x_{ij}\right)\right\}$, $i = 1, \dots, m$, and $\tau_r^{+\prime} = Max\left\{0, \frac{1}{g_r^+}\left(y_{ro} - b\sum_{j\in J}y_{rj}\right)\right\}$, $r = 1, \dots, s$. Then, $(\lambda', \tau^{-\prime}, \tau^{+\prime})$ is a feasible solution for models (6) and (7). □

*Property 2:* Monotonicity [27, 34]

**Theorem 4.3.2** *Let $DMU_p = (a_i x_{io}, i = 1, \dots, m; by_{ro}, r = 1, \dots, s)$ with $0 < a_i \le 1, i = 1, \dots, m$ and $b_r \ge 1, r = 1, .., s$ be a DMU with reduced inputs and enlarged outputs than $DMU_o$. Then, the optimal value to model (6) for $DMU_p$ is not less than that for $DMU_o$, i.e., $\rho_o^F(g) \le \rho_p^F(g)$.*

***Proof*** Obviously, the assertion is true if $a_i = 1$ for all $i = 1, \dots, m$ and $b_r = 1$ for all $r = 1, \dots, s$. Thus, let there exists at least one $i$ or $r$ such that $a_i < 1$ or $b_r > 1$. Without less of generality, let $a_h < 1, h \in \{1, \dots, m\}$. Then, the optimal value to model (6) for DMU$_p$ is obtained by solving the following model:

$$\begin{aligned} \rho_p^F(g) \coloneqq Min\ & \frac{1 + \frac{1}{m}\sum_{i=1}^{m} \tau_i^-}{1 - \frac{1}{s}\sum_{r=1}^{s} \tau_r^+} \\ s.t.\ & \sum_{j \in J} \lambda_j (a_i x_{ij}) \le (a_i x_{io}) + \tau_i^- g_i^-, \quad i = 1, \dots, m, \\ & \sum_{j \in J} \lambda_j (b_r y_{rj}) \ge (b_r y_{ro}) - \tau_r^+ g_r^+, \quad r = 1, \dots, s, \\ & L \le \sum_{j \in J} \lambda_j \le U, \\ & \lambda_j \ge 0, \tau_i^- \ge 0, \tau_r^+ \ge 0, \forall j, i, r. \end{aligned} \tag{16}$$

Given an optimal solution $(\lambda_j^*, \tau_i^{-*}, \tau_r^{+*}, \forall j, i, r)$ to model (16), we have

$$\sum_{j \in J} \lambda_j^* x_{ij} \le x_{io} + \frac{\tau_i^{-*}}{a_i} g_i^-, \quad i = 1, \dots, m,$$

$$\sum_{j \in J} \lambda_j^* y_{rj} \ge y_{ro} - \frac{\tau_r^{+*}}{b_r} g_r^+, \quad r = 1, \dots, s,$$

$$L \le \sum_{j \in J} \lambda_j^* \le U.$$

This indicates that $\left(\lambda_j^*, \frac{\tau_i^{-*}}{a_i}, \frac{\tau_r^{+*}}{b_r}, \forall j, i, r\right)$ is a feasible solution to model (6) for DMU$_o$ with the objective value $\frac{1 + \frac{1}{m}\sum_{i=1}^{m} \frac{\tau_i^{-*}}{a_i}}{1 - \frac{1}{s}\sum_{r=1}^{s} \frac{\tau_r^{+*}}{b_r}}$. Therefore, we have

$$\rho_o^F(g) \le \frac{1 + \frac{1}{m}\sum_{i=1}^{m} \frac{\tau_i^{-*}}{a_i}}{1 - \frac{1}{s}\sum_{r=1}^{s} \frac{\tau_r^{+*}}{b_r}} \le \rho_p^F(g),$$

which completes the proof. □

Proposition 2 in [34] and Theorem 2 in [27] are special cases of this theorem. Note that a similar result holds for $\varphi_o^L(g)$.

*Property 3:* Unit invariance [32, 35]

This property means that the SE scores do not depend on the units of measurement of data.

**Theorem 4.3.3** *Let the direction vector $g$ be such that $g_i^-$, $i = 1, \dots, m$, and $g_r^+$, $r = 1, \dots, s$, have respectively the same unit of the $i$th input and $r$th output. Then, the $\rho_o^F(g)$ and $\rho_o^L(g)$ are all unit invariant.*

***Proof*** Let $(x_j', y_j')$, $j = 1, \dots, n$ be a transformed input-output vector of DMU$_j$ given by:

$$x_{ij}' \coloneqq c_i x_{ij}, y_{rj}' \coloneqq d_r y_{rj}; i = 1, .., m, r = 1, \dots, s, \tag{17}$$

where $c_i, i = 1, \dots, m$ and $d_r, r = 1, \dots, s$ are any collection of positive constants. According to the assumption of theorem, we have

$$g_i^{-\prime} = c_i g_i^-, g_r^{+\prime} = d_r g_r^+; i = 1, \dots, m, r = 1, \dots, s, \tag{18}$$

where $g'$ is the new direction corresponding to the transformed data. By substituting (20) in the constraints of model (6) (or (7)), we have

$$\sum_{j\in J} \lambda_j (c_i x_{ij}) \le c_i x_{io} + \tau_i^- (c_i g_i^-), \qquad i = 1, \dots, m,$$

$$\sum_{j\in J} \lambda_j (d_r y_{rj}) \ge d_r y_{ro} - \tau_r^+ (d_r g_r^+), \quad r = 1, \dots, s.$$

Therefore, the proof is complete. □

For example, the direction vectors (10), (11) and (14) fulfill this property.

*Property 4:* Translation invariance [32, 35-37]

This property indicates that translating the input and output data has no influence on the SE scores. Thus, it is particularly important from computational viewpoint, since it removes the need for data positivity assumption and enables our models to appropriately deal with negative data, accordingly. The following theorem shows that our proposed models can be translation invariant for suitable choices of direction vector.

**Theorem 4.3.4** *The $\rho_o^F(g)$ and $\rho_o^L(g)$ are all translation invariant as long as the specified direction vector is not affected by translation of data.*

***Proof*** Let $(x_j', y_j')$, $j = 1, \dots, n$, be a transformed input-output vector of DMU$_j$ given by:

$$x_{ij}' \coloneqq x_{ij} + a_i, y_{rj}' \coloneqq y_{rj} + b_r; i = 1, .., m, r = 1, \dots, s, \tag{19}$$

where $a_i, i = 1, \dots, m$ and $b_r, r = 1, \dots, s$ are any collection of constants. According to the assumption of theorem, translations of data have no affect on the direction vector so we have

$$g_i^{-\prime} = g_i^-, g_r^{+\prime} = g_r^+; i = 1, \dots, m, r = 1, \dots, s, \tag{20}$$

where $g'$ is the new direction corresponding to the transformed. By substituting (20) in the constraints of model (6) (or model (7)), we have

$$\sum_{j\in J}\lambda_j\big(x_{ij}+a_i\big)\le x_{io}+a_i+\tau_i^- g_i^-,\quad i=1,\dots,m,$$

$$\sum_{j\in J}\lambda_j\big(y_{rj}+b_r\big)\ge y_{ro}+b_r-\tau_r^+ g_r^+,\quad r=1,\dots,s.$$

Since $\sum_{\substack{j=1\\ j\neq o}}^{n}\lambda_j=1$, we eliminate the $a_i$ and $b_r$ on both sides of input and output constraints, and obtain the same constraints as in (6) (or (7)). Therefore, the solution sets are the same and an optimum solution for one program is also optimal for the other. Thus, the proof is complete. □

*Property 5:* Positive, affine transformation invariance [38]

Färe & Grosskopf [38] discussed on a general data transformation, referred to positive, affine data transformation, while measuring the DDF-based efficiency in the VRS framework. This property states that a DEA model is invariant to a positive, affine transformation iff the value of the model does not change for the transformed data given be

$$\begin{aligned} x'_{ij} &:= c_i x_{ij}+a_i,\quad c_i,a_i>0, i=1,..,m;\\ y'_{rj} &:= d_r y_{rj}+b_r,\quad d_r,b_r>0, r=1,\dots,s. \end{aligned}\tag{21}$$

Both the Fractional and Linear GDSE models are positive, affine transformation invariant as long as the chosen direction vector simultaneously satisfies the conditions of Theorems 4.4.3 and 4.4.4. For example, consider the direction vector in (14).

**4.4 Hybrid Directional Super-Efficiency Model**

In some practical applications, the simultaneous proportional/non-proportional contraction (expansion) of *all* inputs (outputs) is not possible. In this sense, some inputs (outputs) are subject to change proportionally while other inputs (outputs) are subject to change non-radially. In such cases to evaluate the SE, we have to deal with a hybrid (see [33] for more details) form of the two groups radial and non-radial SE models. To handle such a situation, integrating RDSE and Linear GDSE models, we formulate a Linear hybrid directional super-efficiency (Linear HDSE) model as follows:

$$Min \quad \frac{1}{m}\left(m_1\tau^- + \sum_{i=m_1+1}^{m} \tau_i^-\right) + \frac{1}{s}\left(s_1\tau^+ + \sum_{r=s_1+1}^{s} \tau_r^+\right)$$

$$s.t. \quad \sum_{j\in J} \lambda_j x_{ij} \le x_{io} + \tau^- g_i^-, \qquad i = 1, \dots, m_1,$$

$$\sum_{j\in J} \lambda_j x_{ij} \le x_{io} + \tau_i^- g_i^-, \qquad i = m_1 + 1, \dots, m,$$

$$\sum_{j\in J} \lambda_j y_{rj} \ge y_{ro} - \tau^+ g_r^+, \qquad r = 1, \dots, s_1,$$

$$\sum_{j\in J} \lambda_j y_{rj} \ge y_{ro} - \tau_r^+ g_r^+, \qquad r = s_1 + 1, \dots, s,$$

$$L \le \sum_{j\in J} \lambda_j \le U, \lambda_j \ge 0, \forall j, \tag{22}$$

$$\tau^- \ge 0, \tau_i^- \ge 0, i = m_1 + 1, \dots, m, \qquad (22.a)$$

$$\tau^+ \ge 0, \tau_r^+ \ge 0, r = s_1 + 1, \dots, s. \qquad (22.b)$$

where

- the direction vector $g$ satisfies the condition (9) in the GRS case (or the condition (13) in the VRS case),
- for $i = 1, \dots, m_1$ and $r = 1, \dots, s_1$, the $i$th input and the $r$th output change proportionally, but with different rates, i.e., $\tau^-$ and $\tau^+$; the data corresponding to these inputs and outputs are supposed to be positive.
- for $i = m_1 + 1, \dots, m$, and $r = s_1 + 1, \dots, s$, the $i$th input and the $r$th output change non-proportionally; the data corresponding to these inputs and outputs are non-negative.

It can be easily proved that the Linear HDSE model is always feasible. Given an optimal solution $(\lambda_j^*, j \in J, \tau^{-*}, \tau_i^{-*}, i = m_1 + 1, \dots, m, \tau^{+*}, \tau_r^{+*}, r = s_1 + 1, \dots, s)$ to model (22), we introduce a Linear HDSE index as follows:

$$\psi_o^{LH}(g) := \left[1 + \frac{1}{m}\left(m_1\tau^{-*} + \sum_{i=m_1+1}^{m} \tau_i^{-*}\right)\right] \times \left[1 - \frac{1}{s}\left(s_1\tau^{+*} + \sum_{r=s_1+1}^{s} \tau_r^{+*}\right)\right]^{-1}.$$

Obviously, the relationship $\rho_o^L(g) \le \psi_o^{LH}(g)$ holds between Linear GDSE and Linear HDSE models.

**4.5 Input-Oriented Directional Super-Efficiency Models**

There are some practical applications, in which the inputs can change and the outputs remain unchanged. To deal with such an application, by setting $g^+ = 0$ in RDSE model in (4), we develop a *radial input-oriented directional super-efficiency* model as follows:

$$\begin{aligned} \gamma_o^I(g) \coloneqq Min \quad & 1 + \tau \\ s.t. \quad & \sum_{j \in J} \lambda_j x_{ij} \leq x_{io} + \tau g_i^-, \quad i = 1, \dots, m, \\ & \sum_{j \in J} \lambda_j y_{rj} \geq y_{ro}, \qquad r = 1, \dots, s, \\ & L \leq \sum_{j \in J} \lambda_j \leq U, \\ & \lambda_j \geq 0, \forall j, \tau{:}\, free\ in\ sign. \end{aligned} \tag{23}$$

Thought this model may be infeasible (e.g., in the VRS case), we provide a necessary and sufficient condition for its feasibility in the special case of $L = 0$, $U = \infty$, i.e., the CRS case.

**Theorem 4.5.1** (*Necessary and Sufficient condition for feasibility*) *The CRS case of model (23) is feasible iff* $Q_o = \emptyset$, *where* $Q_o = \left\{ r \,\middle|\, \sum_{\substack{j=1 \\ j \neq o}}^{n} y_{rj} = 0, y_{ro} > 0 \right\} \subseteq \{1, \dots, s\}$.

***Proof*** (*Sufficiency*) Suppose that $Q_o = \emptyset$. Since $y_o \neq \emptyset$, there exists some index $r$, $r \in \{1, \dots, s\}$, such that $y_{ro} \neq \emptyset$. Therefore, for such $r$, because $Q_o = \emptyset$, there exists a DMU$_j$ such that $y_{rj} \neq \emptyset$. Hence, if $a = \underset{r}{\text{Min}} \left\{ \frac{1}{(n-1) y_{ro}} \sum_{\substack{j=1 \\ j \neq o}}^{n} y_{rj} \,|\, y_{ro} > 0 \right\}$, then $a > 0$. Let $\tau' = \underset{i}{\text{Max}} \left\{ \frac{1}{g_i^-} \left( -x_{io} + \frac{1}{(n-1)a} \sum_{\substack{j=1 \\ j \neq o}}^{n} x_{ij} \right) \right\}$, and $\lambda_j' = \frac{1}{(n-1)a}$, $j = 1, \dots, n$. Obviously, $(\lambda', \tau')$ is a feasible solution to model (23) for the case of CRS. The necessity of the theorem is clear. □

The proposition stated in [39] is a special case of this theorem.

In order to develop a *non-radial input-oriented directional super-efficiency* model considering all the input non-zero slacks, we derive the following model by setting $g^+ = 0$ in Linear NDSE model:

$$\rho_o^I(g) := Min \quad 1 + \sum_{i=1}^{m} \tau_i^-$$
$$s.t. \quad \sum_{j\in J} \lambda_j x_{ij} \leq x_{io} + \tau g_i^-, \quad i = 1, \dots, m,$$
$$\sum_{j\in J} \lambda_j y_{rj} \geq y_{ro}, \qquad r = 1, \dots, s, \tag{24}$$
$$L \leq \sum_{j\in J} \lambda_j \leq U,$$
$$\lambda_j \geq 0, \tau_i^- \geq 0, \forall j, i, \tau: free\ in\ sign.$$

Unlike model (23), increment in inputs may be different for each input in this model. Therefore, the SE score obtained by model (24) is equal or less than that measured by model (23). Clearly, $\rho_o^I(g)$ is always $\geq 1$, and it can be easily verified that only for the extreme-efficient DMUs this index is greater than unity.

Under the CRS assumption, the same condition stated for the feasibility of model (23) ensures the feasibility of model (24).

**Theorem 4.5.2** *(Necessary and Sufficient condition for feasibility) The CRS case of model (24) is feasible iff* $Q_o = \emptyset$.

***Proof*** The feasible solution constructed in the proof of Theorem 4.5.1, is also a feasible one for the CRS version of model (24). Therefore, the proof is the same as proof of Theorem 4.5.1. □

Proposition 3 in [43] is a special case of this theorem.

## 5 Deriving the Conventional Super-Efficiency Models from Our Proposed Models

In this section, we prove that the conventional SE models can be derived from the directional models developed in the preceding section, by assigning specific direction vectors to the latter.

Table 2 shows model (4) includes AP, MAJ, M-MAJ, R-MAJ, Ray and M-Ray models as special cases. Note that in order to derive R-MAJ, Ray and M-Ray models from model (4), the constraints $y_{ro} - \tau g_r^+ \geq 0, r = 1, \dots, s$ are removed.

Table 3 shows that Fractional GDSE model includes Super-SBM-C, Super-SBM-C (I–II), Super-SBM-V and Super-ERM models as special cases. In deriving Super-SBM-I (-V), Super-SBM-C (-V) models, we use the transformations $p_i^- = \tau_i^- g_i^-$ for all $i = 1, \dots, m$ and $q_r^+ = \tau_r^+ g_r^+$ for all $r = 1, \dots, s$. Moreover, in order to derive Super-ERM, we use the transformations $\theta_i = 1 + \tau_i^-$ for all $i = 1, \dots, m$ and $\varphi_r = 1 - \tau_r^+$ for all $r = 1, \dots, s$. As can be seen in Table 3, both Super-SBM-C and Super-ERM models use the same direction vector (10), indicating that these models are equivalent.

**Table 2** The models derived from RDSE model

| Model | Returns to scale | Direction vector |
|---|---|---|
| **AP** | $L = 0, U = \infty$ | $g_i^- = x_{io}, i = 1, \dots, m;\ g_r^+ = 0, r = 1, \dots, s.$ |
| **MAJ** | $L = 0, U = \infty$ | $g_i^- = \underset{j}{\text{Max}}\{x_{ij}\}, i = 1, \dots, m;\ g_r^+ = 0, r = 1, \dots, s.$ |
| **M-MAJ** | $L = 0, U = \infty$ | $g_i^- = \underset{j}{\text{Max}}\{x_{ij} \mid \text{DMU}_j \text{ is efficient}\}, i = 1, \dots, m;\ g_r^+ = 0, r = 1, \dots, s.$ |
| **R-MAJ** | $L = 0, U = \infty$ | $g_i^- = \underset{j}{\text{Max}}\{x_{ij}\}, i = 1, \dots, m;\ g_r^+ = \underset{j}{\text{Max}}\{y_{rj}\}, r = 1, \dots, s.$ |
| **Ray** | $L = 1, U = 1$ | $g_i^- = x_{io}, i = 1, \dots, m;\ g_r^+ = y_{ro}, r = 1, \dots, s.$ |
| **M-Ray** | $L = 1, U = 1$ | $g_i^- = -ax_{io} - 1, i = 1, \dots, m;\ g_r^+ = by_{ro} + 1, r = 1, \dots, s.$ |

$a$ and $b$ are two user-specified parameters as determined in [6].

**Table 3** The models derived from Fractional GDSE model

| Model | Returns to scale | Direction vector |
|---|---|---|
| **Super-SBM-C** | $L = 0, U = \infty$ | $g_i^- = x_{io}, i = 1, \dots, m;\ g_r^+ = y_{ro}, r = 1, \dots, s.$ |
| **Super-SBM-C (I)** | $L = 0, U = \infty$ | $g_i^- = \underset{j}{\text{Max}}\{x_{ij}\}, i = 1, \dots, m;\ g_r^+ = \underset{j}{\text{Max}}\{y_{rj}\}, r = 1, \dots, s.$ |
| **Super-SBM-C (II)** | $L = 1, U = \infty$ | $g_i^- = \underset{j}{\text{Max}}\{x_{ij}\} - \underset{j}{\text{Min}}\{x_{ij}\}, i = 1, \dots, m;\ g_r^+ = \underset{j}{\text{Max}}\{y_{rj}\} - \underset{j}{\text{Min}}\{y_{rj}\}, r = 1, \dots, s.$ |
| **Super-SBM-I** | $L = 0, U = \infty$ | $g_i^- = x_{io}, i = 1, \dots, m;\ g_r^+ = 0, r = 1, \dots, s.$ |
| **Super-SBM-V** | $L = 1, U = 1$ | $g_i^- = x_{io}, i = 1, \dots, m;\ g_r^+ = y_{ro}, r = 1, \dots, s.$ |
| **Super-SBM-I-V** | $L = 1, U = 1$ | $g_i^- = x_{io}, i = 1, \dots, m;\ g_r^+ = 0, r = 1, \dots, s.$ |
| **Super-ERM** | $L = 0, U = \infty$ | $g_i^- = x_{io}, i = 1, \dots, m;\ g_r^+ = y_{ro}, r = 1, \dots, s.$ |

Table 4 gives the details of how to derive LJK, Super-SBM-I, Super-SBM-I-V, Norm1, Norm1-V, Super-Add (I–IV) models from Linear GDSE model. As can be seen in Table 4, both Norm1 and Super-Add (III) models use a same direction vector, indicating that these models are equivalent.

There are some points to be noted:

- In deriving LJK and Super-Add (I–IV) models, the transformations $p_i^- = \tau_i^- g_i^-$ for all $i = 1, \dots, m$ and $q_r^+ = \tau_r^+ g_r^+$ for all $r = 1, \dots, s$ are made in model (7).
- In deriving Norm1 (-V) model, the transformations $x_i = \frac{x_{io}}{\bar{x}_i} + \frac{\tau_i^-}{m}$ for all $i = 1, \dots, m$, and $y_r = \frac{y_{ro}}{\bar{y}_r} - \frac{\tau_r^+}{s}$ for all $r = 1, \dots, s$ are made where $\bar{x}_i = \underset{j}{\text{Max}}\{x_{ij}\}, i = 1, \dots, m$ and $\bar{y}_r = \underset{j}{\text{Max}}\{y_{rj}\}, r = 1, \dots, s$.

*The models derived from Linear HDSE model:*

To derive the Chen et al.'s [51] model from model (22), we suppose that modification in all inputs and outputs be radial, i.e., $m_1 = m$ and $s_1 = s$. Let $L = U = 1$ and the direction vector (10) is used in model (22). Making the transformation $\theta = 1 + \tau^-$ and $\varphi = 1 + \tau^+$, we then have the model of Chen et al. [51], which is infeasible whenever $P_o \neq \emptyset$ or $Q_o \neq \emptyset$.

**Table 4** The models derived from Linear GDSE model

| Model | Returns to scale | Direction vector |
|---|---|---|
| **LJK** | $L = 0, U = \infty$ | $g_i^- = \underset{j}{\text{Max}}\{x_{ij}\}, i = 1, \ldots, m;\ g_r^+ = 0, r = 1, \ldots, s.$ |
| **Norm1** | $L = 0, U = \infty$ | $g_i^- = \frac{1}{m}\underset{j}{\text{Max}}\{x_{ij}\}, i = 1, \ldots, m;\ g_r^+ = \frac{1}{s}\underset{j}{\text{Max}}\{y_{rj}\}, r = 1, \ldots, s.$ |
| **Norm1-V** | $L = 1, U = 1$ | $g_i^- = \frac{1}{m}\underset{j}{\text{Max}}\{x_{ij}\}, i = 1, \ldots, m;\ g_r^+ = \frac{1}{s}\underset{j}{\text{Max}}\{y_{rj}\}, r = 1, \ldots, s.$ |
| **Super-Add (I)** | $L = 0, U = \infty$ | $g_i^- = \frac{1}{m}, i = 1, \ldots, m;\ g_r^+ = \frac{1}{s}, r = 1, \ldots, s.$ |
| **Super-Add (II)** | $L = 0, U = \infty$ | $g_i^- = \frac{m+s}{m}x_{io}, i = 1, \ldots, m;\ g_r^+ = \frac{m+s}{m}y_{ro}, r = 1, \ldots, s.$ |
| **Super-Add (III)** | $L = 0, U = \infty$ | $g_i^- = \frac{m+s}{m}\underset{j}{\text{Max}}\{x_{ij}\}, i = 1, \ldots, m;\ g_r^+ = \frac{m+s}{s}\underset{j}{\text{Max}}\{y_{rj}\}, r = 1, \ldots, s.$ |
| **Super-Add (IV)** | $L = 0, U = \infty$ | $g_i^- = \frac{m+s}{m}\left(\underset{j}{\text{Max}}\{x_{ij}\} - \underset{j}{\text{Min}}\{x_{ij}\}\right), i = 1, \ldots, m;\ g_r^+ = \frac{m+s}{m}\left(\underset{j}{\text{Max}}\{y_{rj}\} - \underset{j}{\text{Min}}\{y_{rj}\}\right), r = 1, \ldots, s.$ |

Now, we turn to the approaches presented by Cook et al. [50] and Chen and Liang [53]. The essence of these approaches is to produce radial/non-radial input- or output-oriented models, by employing the so-called "*big M*" method in non-oriented models.

*Input-oriented models:* Let $L = U = 1$ and $g = (-g^-, g^+) = \left(-x_o, \frac{y_o}{M}\right)$. By removing the constraint $(22.a)$ from model (22), we then have:

(i) If modification in all inputs and outputs be radial, then model (22) reduces to the input-oriented model of Cook et al. [50] by changing variable $\tau^{+\prime} = \frac{\tau^+}{M}$.

(ii) If modification in all inputs be radial, i.e., $m_1 = m$, and in all output be non-radial, i.e., $s_1 = 0$, then model (22) reduces to the input-oriented model of Chen and Liang [53] by changing variables $\tau_r^{+\prime} = \frac{\tau_r^+}{M}$ for all $r = 1, \ldots, s$.

*Output-oriented models*: Let $L = U = 1$, $g = (-g^-, g^+) = \left(-\frac{x_o}{M}, y_o\right)$. By removing the constraint $(22.b)$ from model (22), we then have:

(i) If modification in all inputs and outputs be radial, then model (22) reduces to the output-oriented model of Cook et al. [50] by changing variable $\tau^{-\prime} = \frac{\tau^-}{M}$.

(ii) If modification in all inputs be non-radial, i.e., $m_1 = 0$, and in all outputs be radial, i.e., $s_1 = s$, then model (22) reduces to the output-oriented model of Chen and Liang [53] by changing variables $\tau_i^{-\prime} = \frac{\tau_i^-}{M}$ for all $i = 1, \ldots, m$.

The models presented in [50, 53] suffer from the following shortcomings:

- These model may encounter infeasibility when $P_o \neq \emptyset$ or $Q_o \neq \emptyset$.
- Since each of these models is input- or output-oriented, none of them introduces an integrated non-oriented SE index.

- Theorem 2 proved in [50] (indicates $\beta^* < 1$) is not generally true. Therefore, the presented score for the input-oriented model is not well defined.

## 6 Illustrative Examples

In this section, three examples are given to provide numerical illustrations of comparisons between our proposed measures, with the direction vector (11), and the existing CRS-based and VRS-based SE measures. Using the first example, a comparison between the radial/non-radial CRS super-efficiency models and the Fractional and Linear NDSE models is performed. The second example makes a comparison between our proposed measures and the Norm1-V and the VRS super-efficiency measures developed by [31, 50, 51, 53]. The third example demonstrates that incorporating DM's preference information may affect the super-efficiency scores.

**Example 6.1** Table 5 shows three DMUs using two inputs to produce two outputs. We consider two different instances of output one, denoted by $\mathrm{O}_1$ and $\underline{\mathrm{O}}_1$, where $\underline{\mathrm{O}}_1 = 10 \times \mathrm{O}_1$. Respectively, $\mathrm{O}_1$ and $\underline{\mathrm{O}}_1$ are corresponding to the output one of the $\mathrm{DMU}_j$, $j = 1,2,3$, and $\underline{\mathrm{DMU}}_j$, $j = 1,2,3$. The super-efficiency scores obtained by the conventional CRS super-efficiency models and Fractional and Linear GDSE models in the CRS case are reported in Table 6.

**Table 5** The data set for Example 1

| | $\mathbf{DMU_1}$ | $\mathbf{DMU_2}$ | $\mathbf{DMU_3}$ |
|---|---|---|---|
| $\mathbf{I_1}$ | 1 | 4 | 8 |
| $\mathbf{I_2}$ | 5 | 2 | 1 |
| $\mathbf{O_1}$ | 1 | 0 | 0 |
| $\underline{\mathbf{O}}_\mathbf{1}$ | 10 | 0 | 0 |
| $\mathbf{O_2}$ | 1 | 1 | 1 |

The results can be interpreted as follows:

(A1). Since the obtained scores by R-MAJ model for all DMUs are greater than one, all of them are extreme-efficient and boundary units, as well. Thus, MAJ and M-MAJ models have the same scores.

(A2). It can be observed from Table 2 that, $Q_1 \neq \emptyset$ and $\underline{Q}_1 \neq \emptyset$. Therefore, AP, MAJ, M-MAJ, LJK, Super-SBM-C, Super-SBM (I) and Super-SBM (II) models in evaluating $\mathrm{DMU}_1$ and $\underline{\mathrm{DMU}}_1$ lead to infeasible problems, as expected from Theorems 3, 4 and Remark 2. However, the Fractional and Linear NDSE models are feasible for all DMUs.

(A3). As can be seen in Table 6, among the conventional CRS super-efficiency models, the R-MAJ, Super-Add (I) and Norm1 models are feasible for all DMUs. However, they have other deficiencies that are mentioned below.

(A4). Although R-MAJ model is feasible for all DMUs, in evaluating $\mathrm{DMU}_2$ and $\mathrm{DMU}_3$, it gives the following false projections:

$$\mathrm{DMU}_2 = (4,2,0,1) \longrightarrow \mathrm{Proj}(\mathrm{DMU}_2) = \left(4.6792, 2.4245, \underline{-0.0849}, 0.9151\right),$$

$$\mathrm{DMU}_3 = (8,1,0,1) \longrightarrow \mathrm{Proj}(\mathrm{DMU}_3) = \left(9.1432, 1.7145, \underline{-0.1429}, 0.8571\right).$$

Clearly, in both of these projections the value of the first output is negative so none of them belongs to the corresponding modified PPS.

(A5). It can be observed that all of the feasible models identified $\mathrm{DMU}_1$ as the most efficient DMU.

(A6). Due to the existence of zero in outputs data, the objective functions of Super-SBM-C and Add-Super (II) are undefined for $\mathrm{DMU}_2$ and $\mathrm{DMU}_3$.

(A7). The AP, MAJ, M-MAJ, R-MAJ, Super-SBM-I and Super-Add (I) models, in contrast with the Norm1 and Fractional/Linear GDSE models, have generated a higher score for $\mathrm{DMU}_2$ than $\mathrm{DMU}_3$.

(A8). The Super-Add (I) in evaluating $\underline{\mathrm{DMU}}_1$ yields a large score that arises from its unit variance property.

(A9). The Norm1 model in evaluating $\mathrm{DMU}_1$ yields a score that is greater than unity whereas for other DMUs its scores are less than unity. In addition, there is a considerable difference between the scores obtained from this model for $\mathrm{DMU}_1$ and $\mathrm{DMU}_3$. The cause of this irrationality is that this model objective is not well defined and cannot be interpreted as a super-efficiency index.

(A10). Note that the scores obtained from the Linear GDSE model, are greater than or equal to the scores obtained from the Fractional GDSE model, i.e., $\rho \leq \rho_L^F$.

**Table 6** The results for Example 1

| | **AP** | **MAJ** | **M-MAJ** | **R-MAJ** | **LJK** | **Super-SBM-I** | **Super-SBM-C** | **Super-Add (I)** |
|---|---|---|---|---|---|---|---|---|
| $\mathbf{DMU_1}$ | Inf. | Inf. | Inf. | 2.0000 | Inf. | Inf. | Inf. | 1.7500 |
| $\underline{\mathbf{DMU}}_1$ | Inf. | Inf. | Inf. | 2.0000 | Inf. | Inf. | Inf. | 10.7500 |
| $\mathbf{DMU_2}$ | 1.3000 | 1.1343 | 1.1343 | 1.0849 | 1.2571 | 1.0703 | Und. | 0.2308 |
| $\mathbf{DMU_3}$ | 2.0000 | 1.2000 | 1.2000 | 1.1429 | 1.2000 | 1.2000 | Und. | 0.5000 |

Continued

| | **Super-Add (I)** | **Super-Add (II)** | **Norm1** | **Fractional NDSE** | **Linear NDSE** |
|---|---|---|---|---|---|
| $\mathbf{DMU_1}$ | 1.7500 | Inf. | 1.3750 | 2.3750 | 2.3750 |
| $\underline{\mathbf{DMU}}_1$ | 10.7500 | Inf. | 1.3750 | 2.3750 | 2.3750 |
| $\mathbf{DMU_2}$ | 0.2308 | Und. | 0.2308 | 1.1286 | 1.1304 |
| $\mathbf{DMU_3}$ | 0.5000 | Und. | 0.2000 | 1.1000 | 1.1000 |

**Example 6.2** We now turn to the VRS super-efficiency models. Table 7 shows three DMUs with two inputs and four outputs. The evaluation results of the VRS versions of our proposed measures and the Norm1-V and the former super-efficiency measures developed by [37, 43, 44, 46] are reported in Table 8.

Some interesting observations emerged from these results:

(B1). It can be observed from Table 7 that, $P_1 \neq \emptyset$, $Q_1 \neq \emptyset$ and $Q_4 \neq \emptyset$. Therefore, the models proposed by [37, 43, 44, 46], are infeasible in evaluating $\mathrm{DMU}_1$ and $\mathrm{DMU}_4$. By contrast, our proposed models and Norm1-V models are feasible for all DMUs.

(B2). Similar shortcoming as (A9) can be exposed here for Norm1-V model.

(B3). As can be seen in the Table 8, the ranking orders obtained from Norm1-V, Fractional and Linear NDSE ($\rho_L^F$) models are the same. In addition, $\mathrm{DMU}_1$ is identified as the most efficient, while $\mathrm{DMU}_3$ is regarded as the least efficient by these three measures.

(B4). We note that the Chen and Liang [53] input-oriented model makes no difference between $\mathrm{DMU}_2$ and $\mathrm{DMU}_3$. Cook et al. [50] output-oriented model agrees that, $\mathrm{DMU}_2$ rank is higher than $\mathrm{DMU}_3$ rank. By contrast, the Chen and Liang [53] output-oriented, Cook et al. [50] input-oriented, Chen et al. [51] and Ray [31] models have yield a greater score for $\mathrm{DMU}_2$ than $\mathrm{DMU}_3$. However, none of Cook et al. [50] and Chen and Liang [53] methods provides an integrated non-oriented super-efficiency score for $\mathrm{DMU}_2$ and $\mathrm{DMU}_3$.

To sum up the examples 1 and 2, among all of the conventional CRS/VRS super-efficiency models only the Norm1 and Norm1-V models are feasible. However, the Norm1 and Norm1-V models, in addition to the shortcomings (A9) and (B2), are not able to suitably incorporate the DM's preference information into consideration, as shown in the next example.

**Table 7** The data set for Example 2

| | $\mathrm{DMU}_1$ | $\mathrm{DMU}_2$ | $\mathrm{DMU}_3$ | $\mathrm{DMU}_4$ |
|---|---|---|---|---|
| $\mathrm{I}_1$ | 0 | 1 | 1 | 1 |
| $\mathrm{I}_2$ | 2 | 9 | 1 | 0 |
| $\mathrm{O}_1$ | 1 | 0 | 0 | 0 |
| $\mathrm{O}_2$ | 0 | 4 | 2 | 0 |
| $\mathrm{O}_3$ | 3 | 2 | 1 | 1 |
| $\mathrm{O}_4$ | 0 | 0 | 0 | 1 |

**Table 8** The results for Example 2

| | Norm1-V | Ray 2008 | Chen et al. | Cook et al. I-O | Cook et al. O-O | Chen/Liang I-O | Chen/Liang O-O | Fractional NDSE | Linear NDSE |
|---|---|---|---|---|---|---|---|---|---|
| $\mathbf{DMU_1}$ | 2.5926 | Inf. | Inf. | Inf. | Inf. | Inf. | Inf. | 2.4923 | 2.4923 |
| $\mathbf{DMU_2}$ | 0.7500 | -0.5000 | 2.0000 | 3.0000 | 3.0000 | 4.5000 | 2.0000 | 1.2308 | 1.2308 |
| $\mathbf{DMU_3}$ | 0.3889 | -0.6364 | 4.5000 | 4.5000 | 2.2857 | 4.5000 | 4.5000 | 1.1077 | 1.1077 |
| $\mathbf{DMU_4}$ | 1.1111 | Inf. | Inf. | Inf. | Inf. | Inf. | Inf. | 1.4074 | 1.4074 |

**Example 6.3** Table 9 shows the data for three DMUs each consuming two inputs to produce one output. By using these data, we show that our proposed models have the capability of incorporating the DM's preference information into super-efficiency assessment.

**Table 9** The data set for Example 3

| | $\mathbf{DMU_1}$ | $\mathbf{DMU_2}$ | $\mathbf{DMU_3}$ |
|---|---|---|---|
| $\mathbf{I_1}$ | 1 | 2 | 5 |
| $\mathbf{I_2}$ | 6 | 3 | 2 |
| $\mathbf{O_1}$ | 1 | 1 | 1 |

Suppose that the given weights by the DM corresponding to $\mathrm{I}_1$, $\mathrm{I}_2$ and O respectively are 1, 7 and 1. As explained in the preceding section, to take these information into account, the direction vector (here we use the direction (11)) used in the Fractional NDSE model should be modified as follows:

$$g = (-5,6,1) \Rightarrow g' = \left(-5, \frac{6}{7}, 1\right).$$

The evaluating results by the Norm1 and the CRS version of the Fractional NDSE models, before and after incorporating the DM's preference information are shown in Tables 10 and 11, respectively.

**Table 10** The results for Example 2 before incorporating the DM's preference information

| | **Norm1** | **Fractional NDSE** |
|---|---|---|
| $\mathbf{DMU_1}$ | 0.1000 | 1.1000 |
| $\mathbf{DMU_2}$ | 0.1667 | 1.1667 |
| $\mathbf{DMU_3}$ | 0.0833 | 1.0833 |

**Table 11** The results for Example 2 after incorporating the DM's preference information

| | **Norm1** | **Fractional NDSE** |
|---|---|---|
| $\mathbf{DMU_1}$ | 0.1000 | 1.1000 |
| $\mathbf{DMU_2}$ | 0.1667 | 1.2000 |
| $\mathbf{DMU_3}$ | 0.0833 | 1.5000 |

As can be seen in Tables 10 and 11, the obtained scores for all the DMUs by Norm1 model before and after incorporating the DM's preference information are identical. However, there are noticeable differences between the obtained results by the Fractional NDSE model in these tables. For instance, the rank of $\mathrm{DMU}_3$ before and after taking account the DM's preference information are 3 and 1, respectively. The cause of this happening is that the Fractional NDSE successfully incorporates the DM's preference information, unlike the Norm1. Therefore, our experiment emphasizes the importance of incorporating the DM's preference information along with the PPS's characteristics and illustrates its effect on the ranking order.

## 7 Concluding Remarks

In the present study, employing the DDF, we have systematically examined the super-efficiency DEA models. Using this distance function, we have developed several radial and non-radial input-/ non- oriented directional super-efficiency models and have analyzed each of them, in more details. Our thorough analysis revealed the pivotal role of the direction vector used in these models. Based on this finding, we have demonstrated that each of the conventional super-efficiency models was formulated, indeed, through a directional super-efficiency model, using a specific direction vector; nevertheless, there was no information about this fact and the effects of the direction vector on the model's characteristics such as feasibility, unit invariance, etc. In addition to the role of the direction vector, we have demonstrated the importance of other factors such as (i) the type of the returns to scale (ii) the way used for data adjustment, i.e., radial or non-radial (iii) the orientation of the model. However, the previous studies, which have been reviewed in this paper, did not examine all of these significant findings. Therefore, none of them succeeded to present a complete super-efficiency model and almost all of the conventional models suffer from some shortcomings, specifically the major infeasibility problem.

Finally, our research yielded two complete, named Fractional NDSE and Linear NDSE, measures that are always feasible under the GRS assumption. The proposed models have many additional desirable properties such as stability, unit invariance, capability for dealing with negative data, incorporation of the DM's preference information into ranking and flexibility in computer programming. The results of the illustrative examples have shown the superiority of these models over all of the existing models in possessing the mentioned prosperities.

Because of always feasibility of the proposed models, they are very useful for measuring the Malmquist productivity index. The development of Malmquist indices based on these models is left for future follow-up research.